%
\magnification=\magstep1   
\input amstex
\UseAMSsymbols
\input pictex
\vsize=23truecm
\NoBlackBoxes
\parindent=18pt
  
   \font\rmk=cmr8    \font\itk=cmti8  \font\ttk=cmtt8

\font\gross=cmbx10 scaled\magstep1 

\def\mod{\operatorname{mod}}

\def\Ext{\operatorname{Ext}}

\def\Im{\operatorname{Im}}
\def\bdim{\operatorname{\bold{dim}}}
  
\def\arr#1#2{\arrow <1.5mm> [0.25,0.75] from #1 to #2}

   

 \def\Rahmenbio#1%
   {$$\vbox{\hrule\hbox%
                  {\vrule%
                       \hskip0.5cm%
                            \vbox{\vskip0.3cm\relax%
                               \hbox{$\displaystyle{#1}$}%
                                  \vskip0.3cm}%
                       \hskip0.5cm%
                  \vrule}%
           \hrule}$$}        

\def\Rahmen#1%
   {\centerline{\vbox{\hrule\hbox%
                  {\vrule%
                       \hskip0.5cm%
                            \vbox{\vskip0.3cm\relax%
                               \hbox{{#1}}%
                                  \vskip0.3cm}%
                       \hskip0.5cm%
                  \vrule}%
           \hrule}}}

\vglue3cm
\centerline{\gross Indecomposable representations of the Kronecker quivers}
                    	\bigskip
\centerline{Claus Michael Ringel}     
	\bigskip\bigskip
\plainfootnote{}
{\rmk 2000 \itk Mathematics Subject Classification. \rmk 
Primary 
        16G20. 
Secondary:
05C05, 
11B39, 
15A22, 
16G60, 
17B67, 
65F50. 
}

{\narrower \rmk Abstract.  Let $\ssize k$ be a field
and $\ssize \Lambda$ the $\ssize n$-Kronecker algebra, this is the path
algebra of the quiver with $\ssize 2$ vertices, a source and a sink, and
$\ssize n$ arrows from the source to the sink. It is well-known that the 
dimension vectors of
the indecomposable $\ssize \Lambda$-modules are the positive roots of
the corresponding Kac-Moody algebra. Thorsten Weist has shown 
that for every positive root
there are even tree modules with this dimension vector and that for every
positive imaginary root there are at least $\ssize n$ tree modules. 
Here, we present a short proof of this result.
The considerations used also
provide a calculation-free proof that all exceptional modules over the
path algebra of a finite quiver are tree modules.
\par}
	\bigskip\bigskip
Let $k$ be a field and $Q$ a finite quiver without oriented
cycle. Let $\Lambda= kQ$ be the path algebra of $Q$. 
The target of the paper is to look for $\Lambda$-modules which are tree modules.
According to Kac [K], the 
dimension vectors of the indecomposable $\Lambda$-modules are the positive
roots of the corresponding Lie algebra: for a real root,
there is a unique indecomposable module, for an imaginary root, there are
infinitely many provided $k$ is an infinite field. 
Unfortunately, no effective procedure is known to construct
at least one indecomposable module for each positive root.  On the other
hand, it seems that for each positive root, there exists even a
tree module (the definition will be recalled below): that the 
indecomposable module corresponding to a real root is 
a tree module, and that for
any imaginary root, there are even several different tree modules 
(see [R3], Problem 9).
Thorsten Weist [W] has shown that this is true for all the Kronecker algebras. 
Here, we present a short proof of his result by determining the 
dimension vectors of the ``cover-thin'' Kronecker modules (Proposition 1.1).
	\medskip
The Kronecker algebras are the path algebras of the Kronecker quivers, the
$n$-Kronecker quiver $Q$ with $n$ arrows looks as follows:
$$
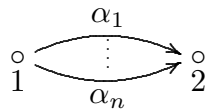
{\beginpicture
\setcoordinatesystem units <1cm,.5cm>
\multiput{$\circ$} at 0 0  2 0 /
\setquadratic
\plot 0.2 0.15  1 0.5  1.8 0.1 /
\plot 0.2 -.15  1 -.5  1.8 -.1 /
\setlinear
\setdots <1mm>
\plot 1 0.4  1 -.4 /
\put{$\alpha_1$} at 1 0,9
\put{$\alpha_n$} at 1 -,9

\put{$1$} at 0 -.5
\put{$2$} at 2 -.5
\setsolid
\arr{1.75 0.15}{1.8 0.1}
\arr{1.75 -.15}{1.8 -.1}
\endpicture}
$$
For $n\ge 2$ we obtain in this way representation-infinite algebras,
for $n\ge 3$ these algebras are wild. The importance of the Kronecker
algebras and their representations is well-known, often they are
considered as the basic data in non-commutative geometry. 
	\medskip
Let $M = (M_a,M_\alpha)_{a,\alpha}$ 
be a finite-dimensional representation of a quiver, thus
$M$ attaches to each vertex $a$ of the quiver a vector space $M_a$
and to each arrow $\alpha$ a linear map $M_\alpha$. 
The sum of the dimension of these vector spaces is called the {\it total
dimension $\dim M$} of $M$. 
In case $M$ is an indecomposable representation with
total dimension $d$, then $M$ is said to be a {\it tree module} provided
that there is a choice of bases for the vector spaces such that the 
corresponding matrix presentations of the linear maps involve altogether 
only $d-1$ non-zero entries (so that the ``coefficient
quiver'' is a tree, see [R2]). 
	\medskip
The root system for the $n$-Kronecker algebra is easy to describe: 
it consists of the vectors $(x,y)\in \Bbb Z^2$ with $x^2+y^2-nxy \le 1.$
The vectors $(x,y)$ with $x^2+y^2-nxy = 1$ are called the real roots, the
other roots the imaginary ones.  
The positive real roots are the dimension vectors of the preprojective
and the preinjective modules. Now, the preprojective and the preinjective
modules are exceptional modules 
(a module over a hereditary algebra
is said to be exceptional if it is indecomposable and has no self-extension)
and exceptional modules over the path algebra of a finite quiver are known to be
tree modules [R2]. Thus, in order to show that every positive root is the
dimension vector of a tree module, we only have to deal with the 
imaginary roots.
	\bigskip
{\bf Theorem.} {\it Let $Q$ be the $n$-Kronecker quiver. 
For any positive imaginary root for $Q$
there are at least $n$ tree modules with this dimension vector.}
	\bigskip
The proof of the theorem will be given in section 2, it will rely on the use
of covering theory. Denote by $\widetilde Q$
the universal covering of $Q$, this is 
the $n$-regular tree with bipartite orientation 
($n$-regular means that every vertex has
precisely $n$ neighbors, the bipartite orientation is characterized by the
property that all vertices are sinks or sources). 
We denote by $\pi\:\mod k\widetilde Q \to \mod kQ$
the push-down functor.
	\medskip
An indecomposable representation of a quiver is said to be {\it thin,} provided
the non-zero vector spaces used are $1$-dimensional. 
If $M$ is a thin 
$k\widetilde Q$-module, then $\pi(M)$ will be said
to be {\it cover-thin.} Similarly, we say that $N$ is {\it cover-exceptional}
provided there is an exceptional $k\widetilde Q$-module $M$ such that
$N = \pi(M)$. 
	\bigskip\bigskip
{\bf 1\. Cover-thin Kronecker-modules.}
	\medskip
{\bf 1.1. Proposition.} {\it Consider $(x,y)\in \Bbb N_0^2$ with $x\le y$.
There exists a cover-thin $kQ$-module $N$ with dimension vector 
$\bdim N = (x,y)$ 
if and only if $0 < y \le (n-1)x+1.$
In this case, 
 there are at least $n$ isomorphism classes of such modules $N$,
unless  $(x,y) = (0,1)$ or $(1,3)$. }

	\bigskip
Proof: Since $k\widetilde Q$ is a tree, the thin indecomposable $k\widetilde Q$-modules
are uniquely determined by the corresponding support, this is just a finite
connected subtree of $k\widetilde Q$. Take a finite connected subtree $T$ of
$k\widetilde Q$ with $|T|$ vertices and let $M(T)$ be the $k\widetilde Q$-module
with support $T$ and $N(T) = \pi(M(T))$
Let $(x,y)$ be the dimension vector of $N(T).$

If $|T| = 1$, then $T$ consists either of 
a sink or a source. The condition $x\le y$ means that $(x,y) = (0,1)$,
that $T$ consists of a sink and that $N(T)$ is the simple projective $kQ$-module.

Now, let $|T|\ge 2.$
Since $T$ is a tree, there is a vertex $a$ in $T$ with
a unique neighbor. In case $a$ is a sink, let $b = a$, otherwise denote by $b$
the unique neighbor of $a$; thus always $b$ is a source. If $b$ is the unique
source in $T$, then $|T| \le n+1$ and $\bdim N(T) = (1,y)$ with $1\le y \le n
= (n-1)+1.$ If $y = n$, then $\pi(M)$ is indecomposable projective and uniquely
determined by its dimension vector, otherwise there are at least $n$ isomorphism
classes of modules of the form $N(T).$

Now assume that $T$ contains at least 2 sources. Removing from $T$ the source $b$
we obtain the disjoint union of a connected tree $T'$ with $|T'| \ge 2$ and
$t\le n-1$ isolated vertices. By induction, we know that $\bdim\pi(M(T')) =
(x',y')$ with $0 < y' \le (n-1)x'+1$ and $(x,y) = (x',y')+(1,t).$ This shows that
$y = y'+t \le (n-1)x'+1+(n-1) = (n-1)x+1.$ (Note that only in case $x' = y'$ and $t=0$,
the pair $(x,y)$ will not satisfy the inequality $x\le y$ we are interested in.)
This shows that the dimension vectors $(x,y)$ of the $kQ$-modules $N(T)$
are as stated. 

Conversely, consider $(x,y)$ with $x\le y$ and 
$0 < y \le (n-1)x+1.$ We try to construct a corresponding $T$. This is clear for
$x \le 1$ and it is easy to see that for $(x,y) = (0,1)$ or $(1,3)$, the
corresponding module $N(T)$ is uniquely determined, wheres for $(1,y)$ with
$1 \le y \le n-1$, there are at least $n$ different isomorphism classes
(for $y = 1$ and for $y = n-1,$ there are precisely $n$ isomorphism classes.  

Thus assume $2 \le x \le y \le (n-1)x+1.$ Write 
$y = \sum\nolimits_{i=1}^x y(i)$ 
with $1 \le y(i) \le n-1$ for $1 \le i \le x-1$ and $1 \le y(x) \le n$
(such a decomposition exists, since $x \le y \le (n-1)x+1$).

Fix some sink $s_1$ of $\widetilde Q$ and take the unique path
$$
 s_1 @<\alpha_1<< t_1 @>\alpha_n>> s_2 @<\alpha_1<< t_2 @>\alpha_n>>\quad \cdots \quad
   @<\alpha_1<< t_{x-1}  @>\alpha_n>> s_x @<\alpha_1<< t_x
$$
starting at $s_1$. For $1\le i \le x,$ we add the arrows $\alpha_j$ (and their
endpoints) starting at $t_i$, with $2 \le j \le y(i).$ 
We see that we obtain in this way a subtree $T$ of $\widetilde Q$, with $x$
sources and $\sum y(i) = y$ sinks, thus $\bdim N(T) = (x,y).$ 

Finally, observe that the module $M(T)$ constructed here for $x \ge 2$ has 
the property that $\Im(\alpha_1)\cap \Im (\alpha_n) \neq 0,$ whereas 
$\Im(\alpha_i)\cap \Im (\alpha_j) = 0$ for $i < j$ and $(i,j)\neq (1,n)$.
Thus, using a permutation of the labels of the arrows, 
the same construction yields $\binom n2$
different isomorphism classes and $\binom n2 \ge n$ for $n \ge 3.$
It remains to consider the case $n = 2.$ Here, $\widetilde Q$ is just a line
and the positive imaginary roots are of the form $(m,m)$ with $m \ge 1.$
Obviously, for every $m\ge 1$, there are precisely two cover-thin $kQ$-modules.
This completes the proof.
	\bigskip
Duality provides in a similar way cover-thin $kQ$-modules with dimension vectors
$(x,y)$ where $0 \le y \le x$ and $0 < x \le (n-1)y+1.$
	\bigskip
It is well-known that the region 
$\Cal F = \{(x,y)\in \Bbb N\mid \frac 1{n-1}x < y \le (n-1)x\}$
is a fundamental domain for the action of the Coxeter transformation on the
set of positive imaginary roots. Note that this region is contained in the
set of dimension vectors of cover-thin $kQ$-modules and that the vectors
$(0,1),\ (1,0),\ (1,3),\ (3,1)$ are real roots. Thus we see:

	\bigskip
{\bf 1.2. Corollary.} {\it For every $(x,y)\in \Cal F$, there are at least $n$
isomorphism classes of cover-thin $kQ$-modules $N$ with $\bdim N = (x,y).$}

	\bigskip
For the benefit of the reader, we provide an illustration for the case $n = 3$:
$$
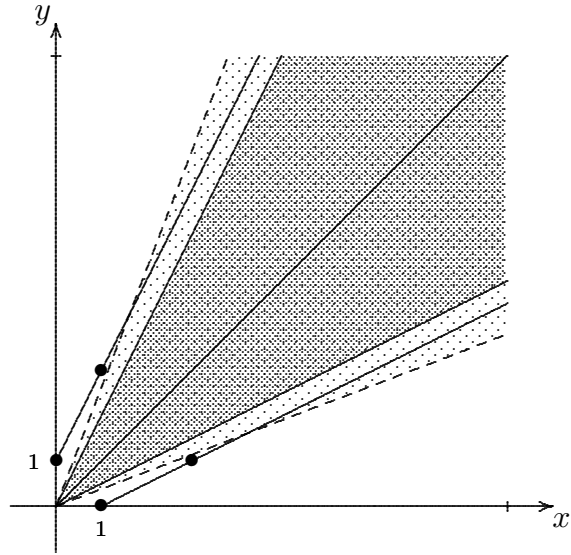
{\beginpicture
\setcoordinatesystem units <.5cm,.5cm>
\arr{-1 0}{11 0}
\arr{0 -1}{0 10.7}
\put{$x$} at 11.2 -.3
\put{$y$} at -.3 10.9
\plot 10 -.1  10 .1 /
\plot -.1 10  .1 10 /
\multiput{$\ssize 1$} at -0.5 1  1 -0.5 /

\plot 0 0  5 10 /
\plot 0 0  10 5 /
\plot 1 0  10 4.5 /
\plot 0 1  4.5 10 /
\plot 0 0  10 10 /
\multiput{$\bullet$} at 1 0  0 1  3 1  1 3    /
\setshadegrid span <0.3mm>
\hshade 0 0 0 <,,,z> 5 2.5 10 <,,z,> 10  5 10 /

\setshadegrid span <.7mm>
\hshade 0 0 0 <,,,z> 3.82 1.4 10 <,,z,> 10  3.82 10 /

\setdashes <1mm>
\plot 0 0  10 3.82 /
\plot 0 0  3.82 10 /
\multiput{$\ssize 1$} at -0.5 1  1 -0.5 /
\endpicture}
$$
The union of the 
shaded areas is the imaginary cone, the dark part being the fundamental domain
$\Cal F$ for the action of the Coxeter transformation on the imaginary cone.
The bullets indicate the dimension vectors $(0,1),\ (1,0),\ (1,3),\ (3,1),$
they are outside of the imaginary cone. There are two lines with slope 2
as well as two lines with slope $\frac 12$: those going through the origin
bound the fundamental region $\Cal F$, the parallel ones bound the region
of the dimension vectors of cover-thin $kQ$-modules.
	\bigskip\bigskip
{\bf 2\. Cover-exceptional $kQ$-modules.}
	\medskip
Thin (indecomposable) modules are exceptional. Thus:
	\medskip
{\bf 2.1. Corollary.} {\it For every positive imaginary root $(x,y)$ 
there are at least $n$
isomorphism classes of cover-exceptional $kQ$-modules  with $\bdim N = (x,y).$}
	\bigskip

{\bf 2.2. Lemma.} {\it Any cover-exceptional module is a tree module.}
	\medskip
Proof: According to [R2], any exceptional module over a hereditary $k$-algebra
is a tree module. But is $M$ is a tree $k\widetilde Q$-module, then
$\pi(M)$ is a tree $kQ$-module. 
	\medskip
The main theorem is a direct consequence of 2.1 and 2.2.
	\bigskip\bigskip
{\bf 3\. Exceptional modules are tree modules.}
	\medskip
The proof of Lemma 2.2 is based on the fact that for $\Lambda$
a finite-dimensional hereditary $k$-algebra, any
exceptional module is a tree module. On the other hand, one can 
use the considerations of section 2
in order to provide a proof of this result which avoids any
calculations.
Indeed, the proof given in [R2] required explicit matrix presentation of 
the preprojective and preinjective Kronecker modules, and, in this way,
was quite technical. Here we show that using induction and the
covering theory for the Kronecker algebras, one can avoid the
matrix calculations.
	\medskip
Using induction on $m,$ we want to show:
	\medskip
{\bf 3.1.} {\it  Let $\Lambda$ be the path algebra of a finite quiver
and $m>0$
a natural number. Any exceptional $\Lambda$-module of length $m$ is a tree module.}
	\medskip 
In the case $m=1$ nothing has to be shown. Thus let us deal with the
induction step, thus let $m > 1.$ 

First, consider the case where $\Lambda$ is the $n$-Kronecker algebra for some
$n\ge 1.$ In the case $n=1$, only one module $N$ has to be considered: it has
length 2 and obviously is a tree module. Thus, assume that $n\ge 2.$ 
The exceptional $\Lambda$-modules are the preprojective modules
$P_0,P_1,P_2,\dots$ 
and the preinjective modules $Q_0,Q_1.Q_2,\dots.$ These modules are
of the form $\pi(M)$ with $M$ an indecomposable representation of
$k\widetilde Q$ (the corresponding dimension vectors in the case $n = 3$
have been displayed in [FR]). Take such a module $M$. Of course, we can assume
that $M$ is not simple. It is easy to see that there is an exact sequence
$$
  0 \to M' \to M \to M'' \to 0
$$
with $\dim\Ext^1(M'',M') = 1$ such that one of the representations 
$M', M''$ is simple and the supports of $M',M''$ are disjoint
(on the level of dimension vectors, we deal with
a subtree $T$ of $\widetilde Q$ and a vertex $a$ in $T$  
having only one neighbor, say $b$; the
coefficients of the dimension vector of $M$ both at $a$ and $b$ are 
equal to $1$. In case $a$ is a sink, $S(a)$ embeds into $M$ say with
image $M'$, then $M'' = M/M'.$ In case $a$ is a source, there is a
surjective map $M \to S(a) = M'',$ in this case $M'$ is chosen as its
kernel. 

Clearly, with $M$ also $M'$ and $M''$ are exceptional modules. Thus,
by induction both are tree modules, and therefore also $M$ is a 
tree module.

Now assume that we are dealing with an exceptional module $M$ of dimension
$M$ such that the support of $M$ has at least three vertices. 
Schofield induction (see [CB] or also [R1]) asserts that there is an exact
sequence
$$
 0 \to X^a \to M \to Y^b \to 0
$$
where $X,Y$ are pairwise orthogonal exceptional modules and the pair
$(a,b)$ is the dimension vector of a sincere preprojective or preinjective
representation $Z$
of an $e$-Kronecker module, with $e = \dim\Ext^1(Y,X).$
Since $a > 0, b> 0$, it follows that $\dim X < m,$ and $\dim Y < m$.
Since the support of $M$ has at least three vertices, we see that not
both modules $X, Y$ can be simple, thus also $\dim Z = a+b < m.$
By induction, all three modules $X,Y, Z$ are tree modules (here, $X,Y$ 
are $\Lambda$-modules, whereas $Z$ is an $e$-Kronecker module),
but then also $M$ is a tree module, see [R2], section 5. This completes
the proof.

	\bigskip\bigskip
\vfill\eject
{\bf References.}
	\medskip
\item{[C]} W. Crawley-Boevey: 
   Exceptional sequences of representations of quivers, 
  in: Representations of algebras, Proc. Ottawa 1992, 
  eds V. Dlab and H. Lenzing, Canadian Math. Soc. Conf. Proc. 14 
  (Amer. Math. Soc., 1993), 117-124. 

\item{[FR]} Ph. Fahr, C.M. Ringel:
   A partition formula for Fibonacci numbers. 
   Journal of Integer Sequences, Vol. 11 (2008), Article 08.1.4. 
    
\item{[K]} V. Kac: Infinite root systems, representations of graphs 
  and invariant theory, Inventiones Math. 56 (1980), 57-92.
\item{[R1]} C.M. Ringel:
   The braid group action on the set of exceptional sequences 
  of a hereditary algebra. 
  In: Abelian Group Theory and Related Topics. Contemp. Math. 171 (1994), 339-352. 
\item{[R2]} C.M. Ringel: Exceptional modules are tree modules. 
   Lin. Alg. Appl. 275-276 (1998) 471-493.
\item{[R3]} C.M. Ringel: 
   Combinatorial Representation Theory: History and Future. 
   In: Representations of Algebras. Vol. I (ed. D.Happel, Y.B.Zhang). BNU Press. 122-144. 
\item{[W]} Th. Weist: Tree modules for the generalized Kronecker quiver. 
   Preprint. \newline ArXiv:0901.1780 [math.RT].

\bigskip\bigskip

{\rmk Fakult\"at f\"ur Mathematik, Universit\"at Bielefeld \par
POBox 100\,131, \ D-33\,501 Bielefeld, Germany \par
e-mail: \ttk ringel\@math.uni-bielefeld.de \par}

\bye